\documentclass[10pt]{article}
 \usepackage{benstyle}

\def \wh {\widehat}
\def \wt {\widetilde}

\def\Frac#1#2{\mbox{\large${\textstyle \frac{#1}{#2}}$}}
\newcommand {\lb} {\label}

\def \ra {{\quad\Rightarrow\quad}}

\title{A stability barrier for reconstructions from Fourier samples}
\author{Ben Adcock \\ Department of Mathematics \\ Purdue University  \\
150 N. University Street \\ West Lafayette, IN 47907 \\ USA \\  \hspace{20pc} \\ \vspace{-2.25pc} \\ \and Anders C. Hansen \\ DAMTP, Centre for Mathematical Sciences \\ University of Cambridge \\ Wilberforce Rd, Cambridge CB3 0WA \\ United Kingdom \and Alexei Shadrin \\ DAMTP, Centre for Mathematical Sciences \\ University of Cambridge \\ Wilberforce Rd, Cambridge CB3 0WA \\ United Kingdom}

\begin{document}
\maketitle

\begin{abstract}
We prove that any stable method for resolving the Gibbs phenomenon---that is, recovering high-order accuracy from the first $m$ Fourier coefficients of an analytic and nonperiodic function---can converge at best root-exponentially fast in $m$.  Any method with faster convergence must also be unstable, and in particular, exponential convergence implies exponential ill-conditioning.  This result is analogous to a recent theorem of Platte, Trefethen \& Kuijlaars concerning recovery from pointwise function values on an equispaced $m$-grid.  The main step in our proof is an estimate for the maximal behaviour of a polynomial of degree $n$ with bounded $m$-term Fourier series, which is related to a conjecture of Hrycak \& Gr\"ochenig.  In the second part of the paper we discuss the implications of our main theorem to polynomial-based interpolation and least-squares approaches for overcoming the Gibbs phenomenon.  Finally, we consider the use of so-called Fourier extensions as an attractive alternative for this problem.  We present numerical results demonstrating rapid convergence in a stable manner.
\end{abstract}

\section{Introduction}
The Fourier series
\bes{
\cF_m f(x) = \frac{1}{\sqrt{2}} \sum_{|j| \leq m} \hat{f}_j \E^{\I j \pi x},\qquad \hat{f}_j = \frac{1}{\sqrt{2}} \int^{1}_{-1} f(x) \E^{-\I j \pi x} \D x,
}
of an analytic and periodic function $f : [-1,1] \rightarrow \bbR$ converges geometrically fast in the truncation parameter $m$.  However, such rapid convergence is destroyed once periodicity is no longer present.  In this case, the series $\cF_m f(x)$ converges only linearly in compact subsets of $(-1,1)$, and there is no uniform convergence on $[-1,1]$.  Near the endpoints $x = \pm 1$ one witnesses the well-known Gibbs phenomenon \cite{GottGibbsRev,Tadmor1}.

Although the Gibbs phenomenon has a long history \cite{HewittGibbs}, and is completely understood mathematically \cite{GottGibbsRev,Tadmor1}, it is a significant hurdle in some applications of Fourier series.  Indeed, it is a testament to the importance of the Gibbs phenomenon that the question of its resolution---that is, restoring high-order convergence from only the first $m$ Fourier coefficients of a function---remains an area of active inquiry.  

Whilst there are many existing methods for this problem (see Section \ref{s:methods} for a review), no one stands out as being inherently superior.  Existing methods with apparent geometric convergence also appear to suffer from some sort of ill-conditioning, and although there are other methods that resolve the Gibbs phenomenon in a stable manner, these typically result in slower orders of convergence.  The purpose of this paper is to explain these observations.  

Specifically, we show that any exponentially-convergent method must also be exponentially ill-conditioned, and in general, if a method has a convergence rate of $\rho^{-m^{\tau}}$ for some $\rho > 1$ and $\tau \in (\frac12,1]$, then it must also possess ill-conditioning of order $\rho^{m^{2 \tau-1}}$.  This result implies the following fundamental stability barrier: the best possible convergence rate for a stable method for resolving the Gibbs phenomenon is root-exponential in $m$.

A theorem of this type is not new.  Our main result is a direct analogue of a theorem of Platte, Trefethen \& Kuijlaars for the problem of overcoming the Runge phenomenon in equispaced polynomial interpolation.  In \cite{TrefPlatteIllCond} it was proved that any method for recovering high accuracy from the pointwise values of a function on an equispaced grid must also exhibit the aforementioned instability behaviour.  The problem we consider in this paper, namely recovering high accuracy from the first $2m+1$ Fourier coefficients of $f$, can be considered a \textit{continuous} analogue of this problem.  Indeed, the problem of recovery from equispaced pointwise values is equivalent to that of recovery from \textit{discrete} Fourier coefficients.

In proving our main result we follow a similar argument to that of \cite{TrefPlatteIllCond}.  The key step therein is the use of an estimate of Coppersmith \& Rivlin \cite{copprivlinpolygrowth} concerning the maximal behaviour of an algebraic polynomial of degree $n$ which is bounded on an equispaced grid of $m$ points.  Our main theoretical contribution is an analogous result for Fourier series: namely, we estimate the maximal behaviour of a polynomial $p$ of degree $n$ with bounded $m$-term Fourier series $\cF_m p$.  In doing so, we provide a partial answer to a conjecture of Hrycak \& Gr\"ochenig \cite{hrycakIPRM} (see Section \ref{s:FCLS} for a discussion).

Our main theorem on the behaviour of polynomials with bounded Fourier sums has an important consequence for a particular method that is sometimes used in practice.  Our result implies that the so-called \textit{inverse polynomial reconstruction method} (IPRM) \cite{shizgalGegen1,shizgalGegen2,pasquettiinverse} is exponentially unstable.  Moreover, although it is possible to stabilize this method via a least squares procedure (henceforth referred to as the \textit{polynomial least squares} method), our main theorem demonstrates that this necessarily decreases the convergence rate to root-exponential. 

Although root-exponential convergence is the best possible permitted, our theorem says nothing about superalgebraic convergence.  Nor does is apply to methods for which convergence occurs only down to some finite tolerance $\epsilon_{\mathrm{tol}}$.  In the final part of this paper we propose the use of so-called \textit{Fourier extensions} \cite{BADHFEResolution,FEStability,BoydFourCont,bruno2003review,brunoFEP,DHFEP} for this problem.  As we discuss, this approach gives rapid convergence---sometimes geometric, but always superalgebraic---but only down to a finite tolerance on the order of machine precision.  We show via example that Fourier extensions typically outperform the aforementioned polynomial least squares method, and conclude that they present an attractive approach to this problem.

\section{The Gibbs phenomenon and its resolution}\label{s:methods}
The Gibbs phenomenon has a long history dating back to Wilbraham in 1848 \cite{wilbraham} (to acknowledge the contribution of Wilbraham, the name Gibbs--Wilbraham phenomenon is also occasionally used).  Forgotten for half a century, this phenomenon was rediscovered by Michelson \cite{michelson}, with the ensuing debate regarding convergence, or lack thereof, between Michelson and Love (carried out in \textit{Nature}) being eventually settled by Gibbs \cite{Gibbsletter1,Gibbsletter2} in 1899, after the arbitration of Poincar{\'e}. The term the Gibbs phenomenon was introduced by B\^{o}cher in 1906 \cite{bocher}.  A detailed review of the Gibbs phenomenon and its history is provided in \cite{HewittGibbs}, with shorter summaries also given in \cite{carlshaw,GottGibbsRev}.

Many methods have been proposed to ameliorate or resolve the Gibbs phenomenon.  Of these, perhaps the earliest to appear were filters and mollifiers.  Here the Gibbs phenomenon is viewed as noise polluting the high-order Fourier coefficients, which can therefore be mitigated by premultiplication with a rapidly decaying function \cite{GottGibbsRev,Tadmor1}.  Unfortunately, standard filters do not lead to high uniform accuracy: they only ensure faster convergence in regions of $[-1,1]$ away from the endpoints $x=\pm 1$.  More recently, Tadmor \& Tanner have developed so-called \textit{adaptive} filters and mollifiers, which lead to greatly improved accuracy \cite{Tadmor1,tadmor2002adaptive,TadmorTanner2005,Tanner2006}.  These can be constructed to obtain geometric convergence in compact subsets of $(-1,1)$ in a stable manner, with typically polynomial accuracy in the vicinity of the endpoints.  For a comprehensive review of this subject, see \cite{Tadmor1}.  Note that this does not contradict the main result of this paper, since the rate of uniform convergence on $[-1,1]$ is not geometric.  As a general principle, geometric convergence away from $x = \pm 1$ can be obtained without ill-conditioning.  

An alternative approach (which can also be viewed as a type of mollifier \cite{Tadmor1}) is the technique of spectral reprojection, introduced and developed by Gottlieb et al.\ \cite{GottGibbs4,GottGibbs3,GottGibbsRev,GottGibbs1,GottGibbsGen}.  Here the slowly convergent Fourier series is reprojected onto a suitable basis; the so-called \textit{Gibbs complementary basis}.  If Gegenbauer polynomials are used for this basis and if the various parameters are carefully tuned, geometric convergence, uniformly on $[-1,1]$, can be restored.  To date, this approach has found application in a number of areas, including image processing \cite{GelbMed2,GelbMed1} and the spectral approximation of PDEs with discontinuous solutions \cite{SMTD}.  

Unfortunately, the original Gegenbauer reconstruction procedure of \cite{GottGibbs1} has been shown to be rather sensitive to the choice of parameters \cite{BoydGegen,GelbTannerGibbs}, with the wrong parameters giving potentially divergent approximations.  To mitigate this effect, a substantially more robust procedure, based on Freud polynomials, was introduced in \cite{GelbTannerGibbs}.  Nonetheless, our main result states that spectral reprojection must either exhibit exponential instability or not be truly exponentially convergent.  We note that, to the best of our knowledge, a stability analysis of spectral reprojection has not yet been carried out.

Spectral reprojection is sometimes referred to a \textit{direct} technique, since it does not require solution of a linear system.  The most obvious \textit{inverse} method is to seek to `interpolate' the first $2m+1$ Fourier coefficients of $f$ with an algebraic polynomial of degree $2m+1$.  This technique, which requires solution of a linear system of equations, is sometimes referred to as the \textit{inverse polynomial reconstruction method} (IPRM) in literature \cite{shizgalGegen1,shizgalGegen2,pasquettiinverse}.  However, interpolation of Fourier coefficients can be seen as a continuous analogue of polynomial interpolation at equispaced nodes.  It should come as little surprise, therefore, that there are substantial issues with both convergence and stability.  In particular, a Runge-type phenomenon is witnessed.  See \cite{BAACHAccRecov,hrycakIPRM}, as well as Section \ref{s:FCLS}, for a discussion.

Since `interpolating' Fourier coefficients may not work, one can also use a lower degree polynomial in combination with a least squares fit (so-called \textit{polynomial least squares}).  This was first discussed in detail in \cite{hrycakIPRM}, and later in \cite{BAACHAccRecov}.  Unfortunately, as we shall prove later,  the degree $n$ can scale at most like $\sqrt{m}$ to ensure stability.  This corresponds to only root-exponential convergence in $m$, consistently with the stability barrier we establish in Section \ref{s:CtsPTK}.  Nonetheless, despite this slower convergence, we remark that this approach does often outperform spectral reprojection in practice.  For a comparison, see \cite{BAACHAccRecov}.

As an alternative to lowering the polynomial degree, one may also try using a higher degree polynomial.  Underdetermined least squares leads to a poor approximation in this case.  However, better accuracy can be restored by using $l^1$-minimization instead.  To the best of our knowledge no analysis currently exists for this approach.  For a related discussion in the case of equispaced function values, see \cite{BoydRunge,TrefPlatteIllCond}.  One may also consider Sobolev norm minimization, such has been considered in the equispaced case in \cite{MinSobolevNorm}.

In \cite{BAACHAccRecov} a general framework was introduced for stable reconstructions in Hilbert spaces.  Given Fourier coefficients, one can reconstruct in any other basis of functions, with one example being the polynomial least squares method discussed above.  An alternative to polynomials involves the use of splines.  As discussed in \cite{BAACHOptimality}, fixed-order splines result in algebraic convergence in a stable manner (see also \cite{WrightFourierSplines}).  One may also consider variable-order splines, but stability becomes an issue.

A different approach to overcoming the Gibbs phenomenon is to smooth the function $f$ via subtraction so as to make it periodic up to a given order, and then compute its Fourier series.  This idea dates back to Krylov and Lanczos, amongst others, and was later studied by Lax and Gottlieb \& Orszag---see \cite{BAthesis,BA3,arnak6} and references therein.  Such smoothing can be carried out implicitly, via extrapolation on the high-order Fourier coefficients; an approach sometimes known as Eckhoff's method in literature \cite{Eckhoff3}.  As shown in \cite{arnak6}, this method converges algebraically fast in $m$.  However, there are also issues with instability \cite{BAthesis}.  A hybrid approach, combining Gegenbauer reconstruction and Eckhoff's method, was also developed in \cite{GottSpliced}.

Alternative methods for overcoming the Gibbs phenomenon arise from sequence extrapolation techniques.  See \cite{BeckermanEpsAlgorithm,BrezSeqAcc}.  In a similar spirit, Driscoll \& Fornberg introduced Pad\'e-based method in \cite{FourPade}.  This approach gives geometric convergence \cite{BeckermanHermitePade}, however in view of our theorem, must also be exponentially unstable.  Such instability was noted in \cite{FourPade}.

There are numerous other methods for resolving the Gibbs phenomenon, and we have not presented a complete list.  The reader is referred to \cite{boydlogsing,GottGibbsRev,Tadmor1} and references therein for further information.  We also mention that many methods designed for the related problem of recovering high accuracy from function values on equispaced grids (i.e.\ overcoming the Runge phenomenon) can potentially be adapted to the Fourier coefficient problem.  See \cite{BoydRunge,TrefPlatteIllCond} for a comprehensive record of such methods.

One such technique that has recently been successfully applied to overcome the Runge phenomenon is that of Fourier extensions.  In the final part of this paper (Section \ref{s:FE}) we consider this approach for overcoming the Gibbs phenomenon.  As we show, the corresponding method can be extremely effective for this problem.

\section{Maximal behaviour of an algebraic polynomial with bounded Fourier series}\label{s:CoppRivFC}
The first step towards the stability theorem in \cite{TrefPlatteIllCond} is an estimate due to Coppersmith \& Rivlin \cite{copprivlinpolygrowth}.  This concerns the behaviour of the quantity
\bes{
A_{n,m} = \sup \left \{ \nm{p}_{\infty} : p \in \bbP_n, \nm{p}_{m,\infty}=1 \right \}, 
}
where $\nm{p}_{\infty} = \sup_{x \in [-1,1]} | p(x) |$, $\nm{p}_{m,\infty} = \max_{j=1,\ldots,m} | p(x_j) |$ and $\{ x_1,\ldots,x_m \}$ is a grid of $m$ equispaced points in $[-1,1]$.  In \cite{copprivlinpolygrowth} it was shown that
\be{
\label{CoppRiv}
(c_1)^{\frac{n^2}{m}} \leq A_{n,m} \leq (c_2)^{\frac{n^2}{m}}.
}
for constants $c_2\geq c_1 >1$.  This estimate determines how many equispaced gridpoints are required to control the behaviour of a polynomial of degree $n$.  Observe that if $m = o(n^2)$ then \R{CoppRiv} implies that there exists a polynomial which is bounded on the grid $\{x_1,\ldots,x_m\}$, but which grows large in between grid points.  On the other hand, $A_{n,m}$ is bounded as $n,m\rightarrow \infty$ if and only if $m = \ord{n^2}$ (that the scaling $m = \ord{n^2}$ is sufficient for boundedness is an older result which dates back to Sch\"onhage \cite{schonhagepolyinterp}).

The problem of the behaviour of $A_{n,m}$ is a classical one in approximation theory (see \cite{TrefPlatteIllCond} for a summary of its history).  However, $A_{n,m}$ involves equispaced grid values.  To prove a stability theorem for overcoming the Gibbs phenomenon, we need to study a related quantity involving Fourier coefficients.  To this end, we define
\be{
\label{Bnm}
B_{n,m} = \sup \left \{ \nm{p}_{2} : p \in \bbP_n, \nm{p}_{m}=1 \right \}, 
}
where $\nm{p}_{2} = \sqrt{\int^{1}_{-1} | p(x) |^2 \D x}$ is the $L^2$-norm on $[-1,1]$ and
\be{
\label{discnormdef}
\nm{p}^2_{m}  = \sum_{|j| \leq m} | \hat{p}_j |^2 =  \nm{\cF_m p}^2_2,
}
is the $l^2$-norm of the first $2m+1$ Fourier coefficients of $p$.

We remark that the quantity $B_{n,m}$ differs from $A_{n,m}$ in two respects.  First, it involves continuous Fourier coefficients, as opposed to discrete Fourier coefficients (i.e. equispaced grid values).  Second, rather than using the uniform norm, we consider the $L^2$ (respectively $l^2$)-norm.  This is natural, since Fourier series satisfy Parseval's relation in the $L^2$/$l^2$-norms.  Indeed, the second equality in \R{discnormdef} follows directly from this relation.  

With this aside, note that 
\bes{
\lim_{m \rightarrow \infty} B_{n,m} = 1,
}
for any fixed $n \in \bbN$.  This follows from strong convergence of the operators $\cF_m \rightarrow \cI$ on $L^2(-1,1)$, and the fact that the space $\bbP_n$ is finite dimensional.  Moreover, much as in the case of $A_{n,m}$, it can be shown that the scaling $m = \ord{n^2}$ is sufficient for boundedness of $B_{n,m}$.  This result was first proved by Hrycak \& Gr\"ochenig \cite[Thm.\ 4.1]{hrycakIPRM} (see also \cite[Lem.\ 3.1]{BAACHAccRecov}).

Our main result in this section shows that $B_{n,m}$ admits a similar lower bound to that found in \R{CoppRiv}.  We have

\thm{
\label{t:Bnmbound}
Let $B_{n,m}$ be as in \R{Bnm}.  Then there exists a constant $c>1$ such that
\bes{
B_{n,m} \geq c^{\frac{n^2}{m}} ,\quad \forall n,m \in \bbN.
}
Specifically, we have the lower bound
\be{
\label{claim}
(B_{n,m})^2 \geq 1+\frac{n}{8 m} + \frac{n}{16 m} d^{\frac{n^2}{m}} ,
}
where $d > \frac{9}{4}$.
}
\prf{
Let $p \in \bbP_n$ be arbitrary.  Integrating by parts, we find that
$$
    \wh p_j 
= \sum_{k=1}^{n} b_k \frac{(-1)^j}{j^{k}}, \qquad  j \in \bbZ \setminus \{0\},
$$
where
$$
    b_k =  -\frac{1}{\sqrt{2}(\I \pi)^{k}}  \left[p^{(k-1)}(1) - p^{(k-1)}(-1)\right].
$$
Therefore, we can write
\be{
\label{polycorres}
   \wh p_j = (-1)^j \wt p(\Frac{1}{j}),\ j \in \bbZ \backslash \{ 0 \}, \qquad  \wt p(t) := \sum_{k=1}^n b_k t^k\,.
   }
Note that $p \in \bbP_n$ is uniquely defined by the values $\hat{p}_0,b_1,\ldots,b_n$.  Therefore, \R{polycorres} defines a one-to-one correspondence between polynomials $p = p(x)$ with $\hat{p}_{0}=0$, say, and polynomials $\tilde p(t)$ satisfying $\tilde p(0) = 0$.  Hence, using Parseval's relation, we have
\eas{
B_{n,m} &\geq \sup \left \{ \nm{p}_2 : p \in \bbP_n, \hat{p}_0=0, \nm{p}_m=1 \right \}  = \sup_{p \in \bbP^*_n} B(n,m,p),
}
where $\bbP^*_n = \{  p \in \bbP_n : p(0) = 0 \}$ and
\bes{
B(n,m,p) = \sqrt{\frac{\sum_{1 \leq | j | < \infty} | p(\frac{1}{j}) |^2}{\sum_{1 \leq | j | \leq m} | p(\frac{1}{j}) |^2} }.
}
To prove the theorem it is sufficient to find a particular $P \in \bbP^*_n$ such that $B(n,m,P)$ admits the bound \R{claim}.

Let $n = 4q+1$ for some $q \in \bbN$.  Consider the polynomial
$$
    P(x) := x T_q^*(x^2) A_{q}(x^2)\,,
$$
where
$$
    A_{q}(x^2) := \prod_{j=1}^{q} (x^2 - \Frac{1}{j^2}),
$$
and $T_q^*(x)$ denotes the Chebyshev polynomial of degree $q$ on the interval $\left[\frac{1}{m^2}, \frac{1}{(q+1)^2}\right]$.  Then, by definition,
$$
     P(\Frac{1}{j}) = 0, \qquad 1 < |j| \le q\,,
$$
and therefore we have
\ea{
     (B(n,m,P))^2 
=  \frac{ \sum_{1 \leq |j| < \infty } |P(\frac{1}{j})|^2 }{ \sum_{1 \leq | j| \leq m} |P(\frac{1}{j})|^2 } 
 =  \frac{ \sum_{q < |j| < \infty } |P(\frac{1}{j})|^2 }{ \sum_{q < |j| \le m} |P(\frac{1}{j})|^2 } 
     = 1 + \frac{ \sum_{m < j<\infty } |P(\frac{1}{j})|^2 } 
             { \sum_{q < j \le m} |P(\frac{1}{j})|^2 } 
     \lb{A}
}
Note that $A_q(\cdot)$ has all its zeros in the interval $[\frac{1}{q^2}, 1]$, hence $|A_q|$ is monotonically decreasing on the interval $[0,\frac{1}{q^2}]$.  Therefore,
$$
     0 < |A_q(\Frac{1}{j_1^2})| 
\le |A_q(\Frac{1}{m^2})| < |A_q(\Frac{1}{j_2^2})|\,, 
     \qquad 
     q < j_1 \le m < j_2\,.
$$
So, putting expression for $P$ into \R{A},  we obtain
\eas{
      (B(n,m,P))^2
 = & 1 + \frac{ \sum_{m < j < \infty } 
       \frac{1}{j^2} |T_q^*(\frac{1}{j^2})|^2  |A_q(\frac{1}{j^2})|^2} 
           { \sum_{q < j \le m} 
       \frac{1}{j_1^2} |T_q^*(\frac{1}{j^2})|^2  |A_q(\frac{1}{j^2})|^2} 
       \\[1ex]
 > & 1 + 
      \frac{ \sum_{m < j  < \infty} 
        \frac{1}{j^2} |T_q^*(\frac{1}{j^2})|^2  |A_q(\frac{1}{m^2})|^2} 
           { \sum_{q < j \le m} 
        \frac{1}{j^2} |T_q^*(\frac{1}{j^2})|^2  |A_q(\frac{1}{m^2})|^2 } 
       \\[1ex]
 = & 1 + 
      \frac{ \sum_{m < j < \infty } 
             \frac{1}{j^2} |T_q^*(\frac{1}{j^2})|^2 } 
           { \sum_{q < j \le m}   
             \frac{1}{j^2} |T_q^*(\frac{1}{j^2})|^2 } \\[1ex]
 := & 1 + \frac{N}{D}
}
Since $|T_q| \le 1$ on the interval $[\Frac{1}{m^2},\Frac{1}{(q+1)^2}]$, 
for the denominator $D$ we have the estimate
$$
    D =  \sum_{q < j \le m}   
             \Frac{1}{j^2} |T_q^*(\Frac{1}{j^2})|^2
      \;\le\; \sum_{q < j < \infty} \Frac{1}{j^2}\,,
$$
so that $D < \frac{1}{q}$.  As for the numerator $N$, we split it into two parts:
\eas{
     N_1 
:= \sum_{m < j < 2m } 
             \frac{1}{j^2} |T_q^*(\Frac{1}{j^2})|^2\,, \qquad N_2
:= \sum_{2m \le j <\infty } 
             \Frac{1}{j^2} |T_q^*(\Frac{1}{j^2})|^2\,.
}
By definition, $T_q^*$ is the Chebyshev polynomial of degree $q$ on the interval $[\frac{1}{m^2},\frac{1}{(q+1)^2}]$, so all of its zeros are located inside that interval. Hence $|T_q^*|$ is decreasing on $[0,\frac{1}{m^2}]$ towards the value 
$|T_q^*(\frac{1}{m^2})| = 1$.  Therefore
\eas{
     N_1
& >  \left |T_q^*\left (\Frac{1}{m^2}\right ) \right |^2 \sum_{m < j < 2m } \frac{1}{j^2}   > \frac{1}{2m}\,, \\
     N_2 
& >  \left |T_q^*\left (\Frac{1}{(2m)^2}\right ) \right |^2  \sum_{2m \le j < \infty } \frac{1}{j^2}  > \Frac{1}{2m} \left |T_q^*\left (\Frac{1}{(2m)^2}\right ) \right |^2\,.
}
i.e.
$$
     N_1 > \frac{1}{2m}\,, \qquad 
     N_2 > \frac{1}{2m}\,\left |T_q^*\left (\Frac{1}{(2m)^2}\right ) \right |^2\,.
$$
Let us evaluate the value $|T_q^*(x^*_0)|$, where $x^*_0 := \Frac{1}{(2m)^2}$.  Take the affine mapping 
$$
    M : \{I^* = [\Frac{1}{m^2},\Frac{1}{(q+1)^2}]\} \to \{I = [-1,1]\}\,,
$$
so that 
$$
   T_q^*(x) = T_q(M(x))  \ra  T_q^*(x^*_0) = T_q(M(x^*_0)) = T_q(x_0)\,,
$$
where $T_q(x)$ is the standard Chebyshev polynomial on $I=[-1,1]$.  The length of the interval $I^*$ is 
less than $\frac{1}{q^2}$, so mapping $M$ onto the interval $I$ with the length $2$ uses the length magnifying factor $\lambda > 2q^2$.   The point $x^*_0 = \frac{1}{(2m)^2}$ lies at the distance $\delta^*_0 = \frac{3}{4m^2}$ from the left endpoint $\frac{1}{m^2}$ of $I^*$, so it will be mapped to the point $x_0 < -1$ given by
$$
    - x_0 = 1 + \delta_0, \qquad 
    \delta_0 > 2q^2 \delta^*_0 = \Frac{3q^2}{2m^2}\,.
$$
For $x = 1 + \delta$, we have
\eas{
      T_q(x) 
& =  \Frac{1}{2}\Big((x + \sqrt{x^2-1})^q + (x - \sqrt{x^2-1})^q\Big) 
>  \Frac{1}{2} (x + \sqrt{x^2-1})^q 
 >  \Frac{1}{2}( 1 + \sqrt{2\delta})^q\,,
}
and since $\sqrt{2\delta_0} = \sqrt{3}\frac{q}{m} > \frac{q}{m}$,
we have
\ea{
    T_q^*(\Frac{1}{(2m)^2}) = T_q^*(x_0^*) = T_q(x_0) > \Frac{1}{2}\Big(1 + \Frac{q}{m}\Big)^q 
    = \Frac{1}{2}\Big[\big(1 + \Frac{q}{m}\big)^{\frac{m}{q}}\Big]^{q^2/m} := \Frac{1}{2} \gamma^{q^2/m}\,.
}
Hence
\be{
   N_2 > \frac{1}{4m}  \gamma^{2q^2/m}, \qquad
   \gamma := (1 + \Frac{1}{r})^r, \qquad   r := \frac{m}{q}\,.
\lb{N_2}
}
Combing these results together, and using the fact that $q \approx \frac{n}{4}$, we obtain
\eas{
    (B(n,m,P) )^2
> & 1 + \frac{N}{D} =  1 + \frac{N_1}{D} + \frac{N_2}{D} \\
 > & 1 + \frac{q}{2m} + \frac{q}{4m} \gamma^{2q^2/m} \\
 > & 1 + \frac{n}{8m} + \frac{n}{16m} \gamma^{n^2/8m}\,.
}
Since $m > n/2$ and $q < n/4$, we have $ r= \frac{m}{q} > 2$ so that a lower bound for $\gamma$ is $\gamma  = (1 + \Frac{1}{r})^r > 9/4$.  This completes the proof.
}

\begin{figure}
\begin{center}
$\begin{array}{ccc}
 \includegraphics[width=4.75cm]{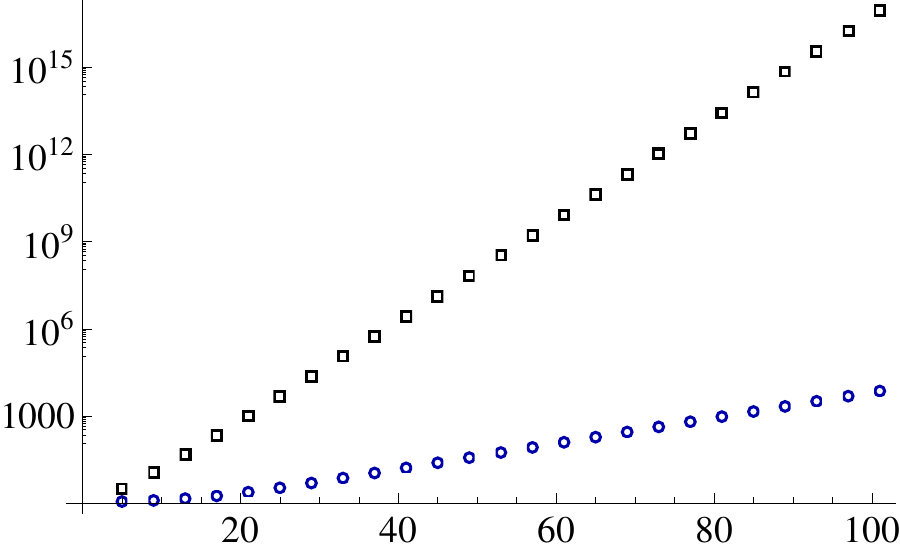}   &   \includegraphics[width=4.75cm]{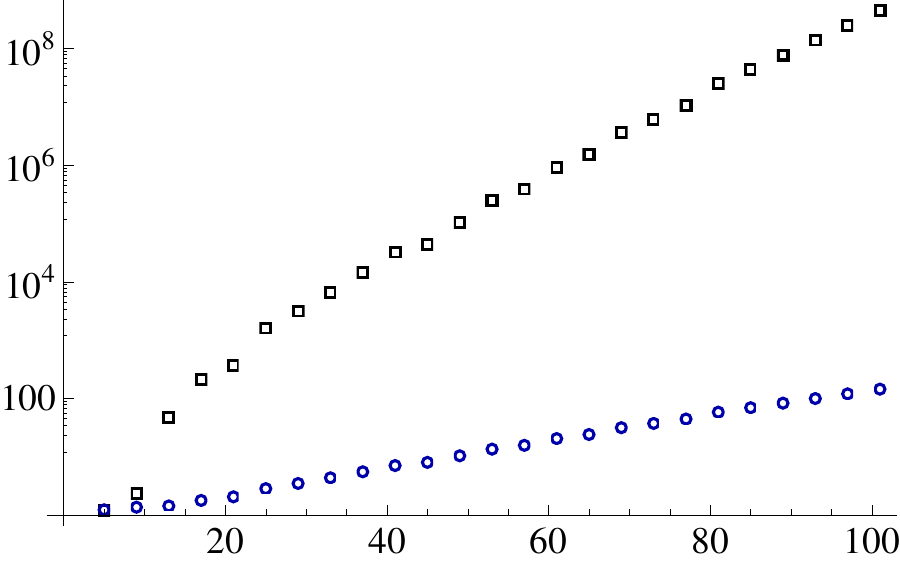} &   \includegraphics[width=4.75cm]{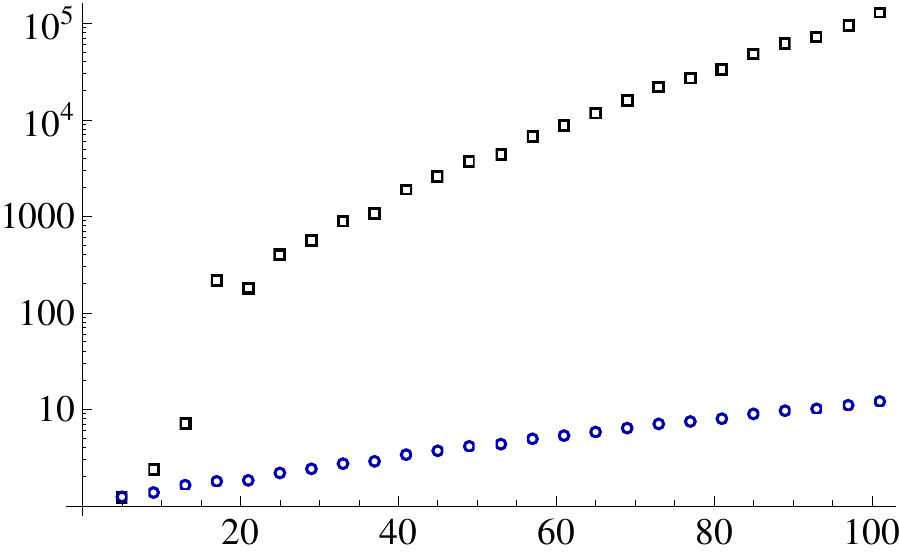} \\  \\\includegraphics[width=4.75cm]{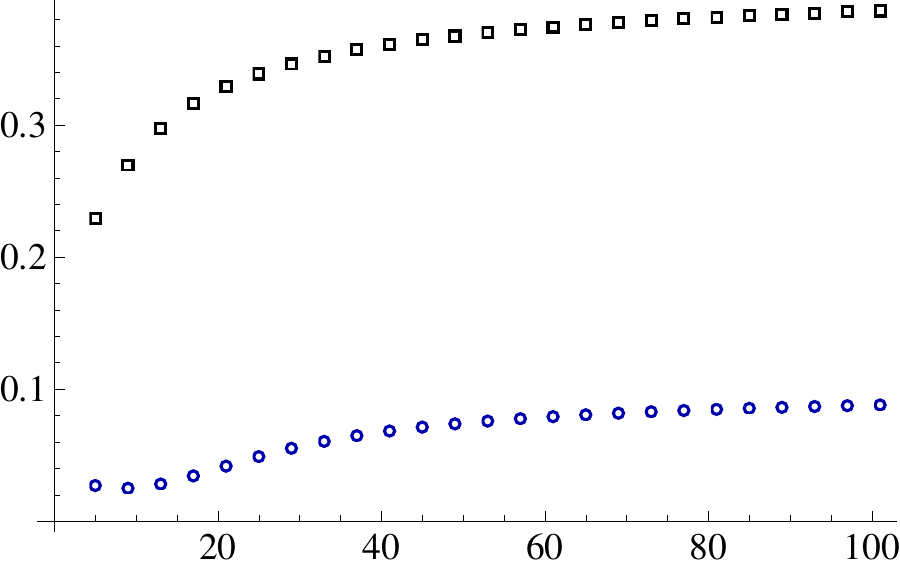}   &   \includegraphics[width=4.75cm]{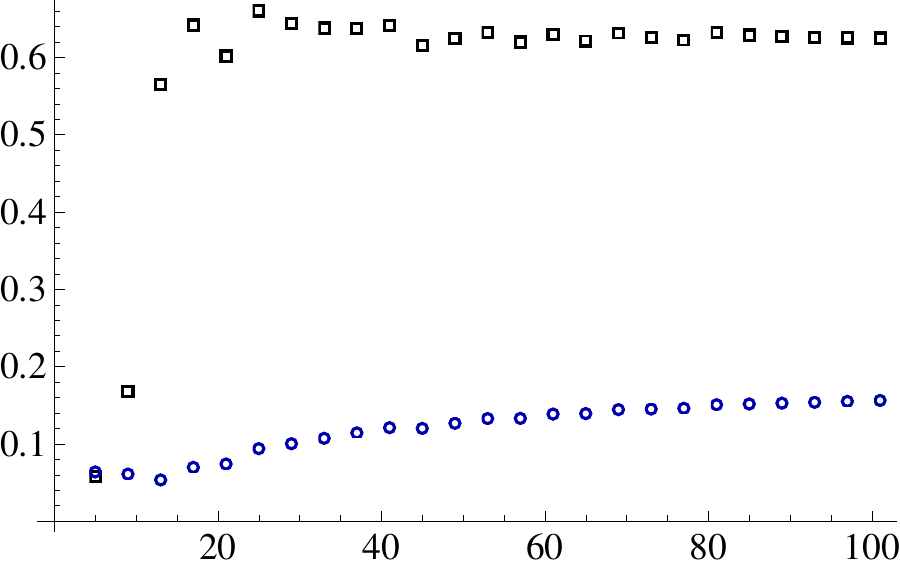} &   \includegraphics[width=4.75cm]{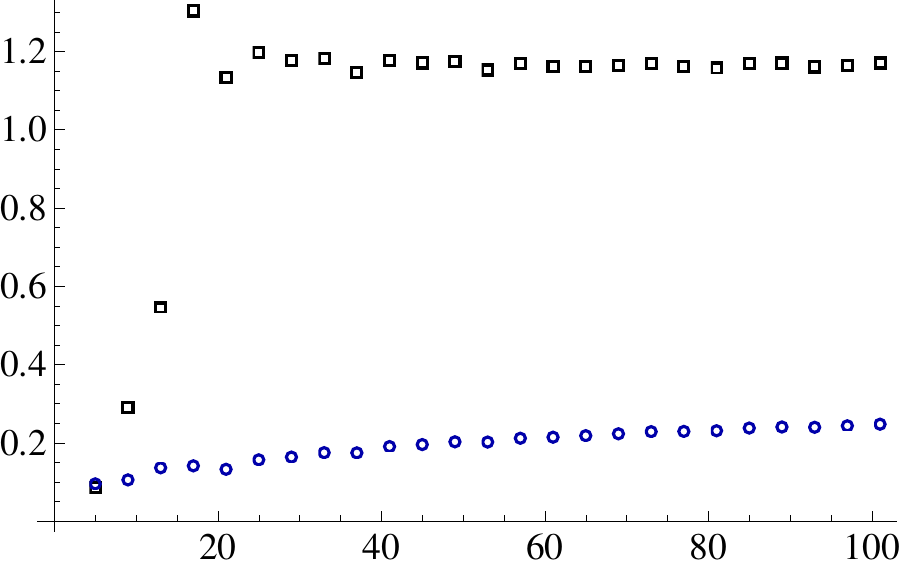} \\ \\
 (\alpha,\beta) = (\frac12,1) &  (\alpha,\beta) = (\frac14,\frac54) &  (\alpha,\beta) = (\frac18,\frac32)
\end{array}$
\caption{\small Top row: the quantity $B_{n,\alpha n^{\beta}}$ (squares) and the lower bound $B^*_{n,\alpha n^{\beta}}$ (circles) against $n$. Bottom row: the scaled quantities $n^{\beta - 2} \log (B_{n,\alpha n^{\beta}} )$ and $n^{\beta - 2} \log (B^*_{n,\alpha n^{\beta}})$.  Computations were carried out in \textit{Mathematica} using additional precision.  Note that the quantity $B_{n,m}$ can be computed, since coincides with the minimum singular value of a particular matrix.  }  \label{f:boundcomparison}
\end{center}
\end{figure}

In Figure \ref{f:boundcomparison} we confirm this theorem by plotting the quantity $B_{n,m}$ and the lower bound
\bes{
B^*_{n,m} = \sqrt{1+\frac{n}{8 m} + \frac{n}{16 m} \left ( \frac{9}{4} \right )^{\frac{n^2}{m}} }.
}
Note that $B^{*}_{n,m}$ not only provides a lower bound, it also appears to correctly predict the behaviour of $B_{n,m}$.  We conjecture that an upper bound of the form
\be{
\label{conjecture}
B_{n,m} \leq c^{\frac{n^2}{m}},\quad \forall n,m \in \bbN,
}
also holds, much as in the case of the quantity $A_{n,m}$.

\section{An impossibility theorem for the resolution of the Gibbs phenomenon}\label{s:CtsPTK}
We now turn our attention to the main result of this paper.  First, we require some notation.  Following \cite{TrefPlatteIllCond}, let $\{ \phi_m \}_{m \in \bbN}$ be a family of mapping $L^2(-1,1) \rightarrow L^2(-1,1)$ such that $\phi_m(f)$ depends only on the values $\{ \hat{f}_j \}_{|j| \leq m}$.  Note that the mappings $\phi_m$ can be both linear or nonlinear.  We define the condition number for $\phi_m$ by
\bes{
\kappa_m = \sup_{f} \lim_{\epsilon \rightarrow 0} \sup_{\substack{g : \\ 0 < \nmu{g}_m \leq \epsilon }} \frac{\| \phi_m(f+g) - \phi_m(f) \|_2}{ \nmu{g}_m }.
}
For a compact set $E \subseteq \bbC$ we shall also let $B(E)$ be the Banach space of functions continuous on $E$ and analytic in its interior, with norm $\nm{f}_E = \sup_{z \in E} | f(z) |$.

Our main result is as follows:

\thm{
\label{t:ctsPTK}
Let a compact set $E \subseteq \bbC$ containing $[-1,1]$ in its interior be fixed, and suppose that $\{ \phi_m \}_{m \in \bbN}$ are such that (i) for any $f$, the quantity $\phi_m(f)$ depends only on the values $\{ \hat{f}_j \}_{|j| \leq m}$, and (ii) for some $M<\infty$, $\sigma > 1$ and $\tau \in (\frac{1}{2},1]$, we have
\be{
\label{E:Phiferr}
\| f - \phi_m(f) \|_2 \leq M \sigma^{-m^{\tau}} \| f \|_{E},\qquad \forall f \in B(E), m \in \bbN.
}
Then the condition numbers $\kappa_m$ satisfy
\bes{
\kappa_m \geq c^{m^{2 \tau-1}},
}
for some $c>1$ and all sufficiently large $m$.
}
\prf{
We proceed as in \cite{TrefPlatteIllCond}.  Without loss of generality we may assume that $E$ is a Bernstein ellipse $E(\rho)$ for some $\rho > 1$.  We can now replace \R{E:Phiferr} by
\be{
\label{step1}
\| f - \phi_m(f) \|_2 \leq \frac{1}{2} \rho^{-\alpha m^{\tau}} \nm{f}_E,\quad m\geq m_0,\quad \forall f \in B(E(\rho)).
}
where $m_0$ is sufficiently large and $\alpha > 0$ is fixed.  Let $p \in \bbP_n$.  An inequality of Bernstein \cite[Lem.\ 1]{TrefPlatteIllCond} implies that
\bes{
\nm{p}_E \leq \rho^{n} \| p \|_{[-1,1]}.
}
Since $p \in \bbP_n$, we have that $\| p \|_{[-1,1]} \leq c n \| p \|_2$ for some constant $c > 0$.  Thus,
\be{
\label{crude}
\| p \|_{E} \leq c n \rho^{n} \| p \|_2.
}
Setting $p = f$ in \R{step1} now gives
\eas{
\| \phi_m(p) \|_2  \geq \| p \|_2 - \tfrac{1}{2} \rho^{-\alpha m^{\tau}} \| p \|_{E}  \geq \left ( 1 - \tfrac{1}{2} c n \rho^{n- \alpha m^{\tau}} \right ) \| p \|_2.
}
Suppose that $n \leq \frac{1}{2} \alpha m^{\tau}$.  Then
\bes{
c n \rho^{n-\alpha m^{\tau}} \leq \tfrac{1}{2} c \alpha m^{\tau} \rho^{-\frac{1}{2} \alpha m^{\tau}} < 1,\quad \forall m \geq m_1,
}
where $m_1$ is sufficiently large.  Set $m^* := \max \{ m_0,m_1 \}$.  If $m \geq m^*$ this now gives
\bes{
\| \phi_m(p) \|_2 \geq \tfrac{1}{2} \| p \|_2,\quad \forall m \geq m^*.
}
Since $\phi_m(0) = 0$ (this follows from \R{E:Phiferr}) we have
\bes{
 \frac{\| \phi_m(\epsilon p) - \phi_n(0) \|_2}{\nm{\epsilon p }_m} = \frac{\| \phi_m(\epsilon p)  \|_2}{\nm{\epsilon p }_m} \geq \frac{1}{2} \frac{\| \epsilon p \|_2}{\nm{ \epsilon p}_m} = \frac{1}{2} \frac{\nm{p}_2}{\nm{p}_m}.
}
Therefore, since $p \in \bbP_n$ is arbitrary, we obtain
\bes{
\kappa_m \geq \tfrac{1}{2} \sup_{\substack{p \in \bbP_{n} \\ p \neq 0}}  \left \{ \frac{\nm{p}_2}{\nm{p}_m}  \right \}= \tfrac12 B(n,m),\quad \forall n \leq \tfrac{1}{2} \alpha m^{\tau},
}
where $B(n,m)$ is given by \R{Bnm}.  Using Theorem \ref{t:Bnmbound} with $n = \frac{1}{2} \alpha m^{\tau}$, this now yields
\bes{
\kappa_m \geq \tfrac{1}{2} c^{\frac{1}{4} \alpha^2 m^{2 \tau-1}},
}
as required.
}

This theorem implies that any method for overcoming the Gibbs phenomenon that results in a convergence rate of $\sigma^{-m^{\tau}}$ for all analytic functions $f \in B(E)$ must also possess ill-conditioning of order $c^{m^{2 \tau-1}}$.  In particular, the best convergence rate that can be achieved with a stable method is root-exponential in the number of Fourier coefficients $m$.  As commented in \cite{TrefPlatteIllCond}, we stress that, despite the use of polynomials in the proof, this theorem is not about polynomials or polynomial-based approximation procedures.  It holds for all (linear or nonlinear) mappings $\phi_m$ satisfying a bound of the form \R{E:Phiferr}.

\section{Fourier coefficient interpolation and least squares}\label{s:FCLS}
As discussed in Section \ref{s:methods}, an obvious way to seek to overcome the Gibbs phenomenon is by interpolating the first $2m+1$ Fourier coefficients of $f$ with a polynomial of degree $2m$.  In other words, we construct an approximation $f_m$ satisfying
\be{
\label{IPRM}
\widehat{f_m}_{j} = \hat{f}_j,\quad |j| \leq m,\qquad f_m \in \bbP_{2m}.
}
It can be shown that such a polynomial exists uniquely for any $m$ \cite{hrycakIPRM}.  However, as we shall see in a moment, this approach is both exponentially unstable and divergent.  A simple modification involves computing an overdetermined least squares with ($n \leq m$):
\be{
\label{LS}
f_{n,m} =  \underset{p \in \bbP_{2n}}{\operatorname{argmin}} \Bigg \{ \sum_{|j| \leq m} | \hat{f}_j - \hat{p}_j |^2 \Bigg \}.
}
Note that $f_m$, as defined by \R{IPRM}, coincides with $f_{n,m}$ when $n=m$.  As mentioned, \R{IPRM} is often referred to as the inverse polynomial reconstruction method (IPRM) \cite{shizgalGegen1,shizgalGegen2}.  The modification \R{LS}, referred to as polynomial least squares, was introduced and analysed in \cite{hrycakIPRM}, and developed further in \cite{BAACHAccRecov} (see also \cite{AdcockHansenSpecData}).

In \cite{BAACHOptimality} it was shown that $f_{n,m}$ satisfies the sharp bound
\bes{
\| f - f_{n,m} \|_2 \leq B_{2n,m} \inf_{p \in \bbP_{2n}} \| f - p \|_2,
}
where $B_{n,m}$ is the constant defined in \R{Bnm} (previous, but non-sharp, estimates were given in \cite{BAACHAccRecov,hrycakIPRM}).   Moreover, it was also shown that the condition number $\kappa = \kappa_{n,m}$ of the mapping $f \mapsto f_{n,m}$ is precisely
\be{
\label{kappa_equality}
\kappa_{n,m} = B_{2n,m},
}
where $B_{n,m}$ is defined by \R{Bnm}.  Hence, Theorem \ref{t:Bnmbound} allows us to explain both the convergence and stability of this approach.  In particular, \R{kappa_equality} and Theorem \ref{t:Bnmbound} gives that
\be{
\label{OurGrochenig}
\kappa_{n,m} \geq c^{\frac{n^2}{m}}.
}
meaning that polynomial least squares is unstable whenever $n$ grows faster than $\sqrt{m}$.  In the particular case $n=m$, this shows that the IPRM method \R{IPRM} is exponentially ill-conditioned.  Note that such exponential growth was previously observed, although not analysed, in \cite{hrycakIPRM,shizgalGegen2}.  

We remark that Hrycak \& Gr\"ochenig have conjectured that an upper bound of the form \R{OurGrochenig} also holds.  In other words, the condition numbers can grow no worse than $c^{\frac{n^2}{m}}$ for some $c>1$.  This is of course equivalent to the conjecture \R{conjecture} concerning the constant $B_{n,m}$.

Regarding the convergence of $f_{n,m}$, it is useful to recall that the error $\inf_{p \in \bbP_{2n}} \| f - p \|$ decays like $\rho^{-2n}$, where $\rho > 1$ is the parameter of the largest Bernstein ellipse within which $f$ is analytic \cite{rivlin1990chebyshev}.  Thus, we have the bound
\be{
\label{expbound}
\| f - f_{n,m} \|_2 \leq c_f B_{n,m} \rho^{-2n},
}
where $c_f > 0$ is some constant depending on $f$ only.  This indicates that $f_{n,m}$ may fail to converge to $f$ if the rate of growth of $B_{n,m}$ exceeds that of $\rho^{2 n}$.  In particular, if $n = \cO(m)$, in which case $B_{n,m}$ is exponentially large, there will always be functions (analytic within only a small ellipse $E(\rho))$ for which the right-hand side of \R{expbound} diverges.  Thus, the polynomial least squares method \R{LS}, and in particular the IPRM \R{IPRM},  may well suffer from a Runge-type phenomenon---i.e.\ divergence of $f_{n,m}$ for some nontrivial family of analytic functions---whenever $n = \ord{m}$.  Although we have no proof of this fact (the upper bound \R{expbound} need not be sharp for an individual function $f$), there is substantial numerical evidence to support this conjecture \cite{hrycakIPRM,shizgalGegen2}.

On the other hand, when $n=\ord{\sqrt{m}}$ it is known that $B_{n,m} = \ord{1}$ \cite{BAACHAccRecov,hrycakIPRM}, and therefore we have stability, as well as guaranteed convergence of the approximation.  Unfortunately, the convergence rate is only root-exponential, and Theorem \ref{t:ctsPTK} demonstrates that it can be no faster.
 
\rem{
Overdetermined least squares is a well-known approach to overcome the Runge phenomenon in equispaced polynomial interpolation \cite{boyd2009divergence,TrefPlatteIllCond}.  In \cite{boyd2009divergence} essentially the same arguments as those given above were presented for this problem, leading to the same conclusions: namely, a Runge-type phenomenon for $n = \ord{m}$, but stability and root-exponential convergence whenever $n = \ord{\sqrt{m}}$.

We note that the use of polynomial least squares for overcoming the Gibbs phenomenon, as opposed to the Runge phenomenon, is far less well known.  However, this approach (with $n = \ord{\sqrt{m}}$) does appear to outperform other more commonly used techniques, such as spectral reprojection, despite its formally slower convergence \cite{BAACHAccRecov}.
}

\section{Fourier extensions for overcoming the Gibbs phenomenon}\label{s:FE}
The principle behind polynomial least squares is to reconstruct a function $f$ in a finite-dimensional subspace in which it is well-approximated, with the space of polynomials of degree $n$ being a natural choice for analytic $f$.  Recently these ideas were developed substantially in several papers by Adcock \& Hansen \cite{BAACHAccRecov,BAACHOptimality}.  Therein a framework, \textit{generalized sampling}, was introduced to stably reconstruct elements of Hilbert spaces in finite-dimensional subspaces from their samples taken with respect to a particular basis or frame.  Polynomial least squares is a specific instance of this general framework, based on sampling with respect to the Fourier basis and reconstructing in the subspace $\bbP_{2n}$.  However, since generalized sampling allows reconstructions in arbitrary spaces, there is no need to choose this particular space.  In this section, we consider the use of an alternative subspace for reconstruction, and demonstrate that this gives substantial improvements over polynomial least squares. 

The particular subspace we shall employ is based on an approximation scheme known as \textit{Fourier extensions}.  Here one seeks to approximate a function using a Fourier series on the extended domain $[-T,T]$, where $T>1$ is fixed.  In other words, we compute an approximation belonging to the space
\be{
\label{Sn_def}
\cS_n = \left \{ \sum_{|j| \leq n} a_j \E^{\I \frac{j \pi}{T} x} : a_j \in \bbC, | j | \leq n \right \}.
}
Fourier extensions have been used in the past to successfully overcome the Runge phenomenon.  In particular, Boyd \cite{BoydFourCont,BoydRunge}, Bruno \cite{bruno2003review} and Bruno et al.\ \cite{brunoFEP} have shown that this approximation can recover analytic functions to extremely high accuracy from equispaced data.  This was confirmed by the analysis of Huybrechs \cite{DHFEP}, and later Adcock et al.\ \cite{FEStability}.

Computing a Fourier extension from equispaced data requires solving an extremely ill-conditioned linear system.  Hence there can be substantial differences between the exact (i.e.\ infinite-precision) and numerical (i.e.\ finite-precision) approximation.  Both were analysed in \cite{FEStability} in the context of the impossibility theorem of \cite{TrefPlatteIllCond}.  Specifically, it was shown that
\begin{enumerate}
\item[(i)] In infinite precision, the Fourier extension computed from equispaced data attains the stability barrier of \cite{TrefPlatteIllCond}.  In other words, it is stable and converges root-exponentially fast in $m$.
\item[(ii)] In finite precision, the corresponding approximation $f_m$ is numerically stable and converges geometrically fast for all functions analytic in a sufficiently large complex region $E$, but only down to a finite accuracy on the order of $\epsilon$, where $\epsilon$ is proportional to the machine precision used.  Specifically, one has a bound of the form
\be{
\label{Fin_Err_Bound}
\| f - f_m \|_2 \leq M \left ( \sigma^{-m} + \epsilon \right ) \| f \|_E.
}
For all other analytic functions, the convergence is at least superalgebraic down to the same accuracy.
\end{enumerate}
These results, in particular (ii), make Fourier extensions effective methods for the equispaced data problem.  Note that there is no contradiction with the impossibility theorem of \cite{TrefPlatteIllCond}, since the bound \R{Fin_Err_Bound} contains the finite term $\epsilon$ (compare with \R{E:Phiferr}).  Fortunately, this term is on the order of machine precision.  Hence it has no substantial effect on the overall approximation whilst simultaneously allowing the stability barrier to be circumvented.

We now consider the use of Fourier extensions for reconstructions from Fourier samples.  We shall consider the following approximation:
\be{
\label{FE}
\tilde{f}_{n,m} = \underset{\phi \in \cS_{n}}{\operatorname{argmin}} \Bigg \{ \sum_{|j| \leq m} | \hat{f}_j - \hat{\phi}_j |^2 \Bigg \}.
}
Note this construction is almost identical to that considered in \cite{FEStability,BoydRunge,brunoFEP} for the case of equispaced data, the only difference being that the discrete least-squares is taken over a sum of Fourier coefficients as opposed to pointwise function values.  Note also that $\tilde{f}_{n,m}$ is similar to the polynomial least-squares approximation \R{LS}, but with the approximation space $\cS_n$ as opposed to $\bbP_{2n}$.

It is possible to adapt the analysis of \cite{FEStability} for equispaced data to the case of \R{FE}.  For the sake of brevity, we shall not do this.    Instead, we shall show by numerical example how effective this approximation can be in practice.  However, we remark that one can show that conclusions (i) and (ii) stated above also hold in this case.  In particular, the numerical approximation obtained in finite precision effectively avoids the stability barrier of Theorem \ref{t:ctsPTK}.

To demonstrate the effectiveness of \R{FE} we shall compare it to the polynomial least squares approximation $f_{n,m}$ given by \R{LS}.  Both methods involve additional parameters: $n$, the polynomial degree, in the case of $f_{n,m}$, and $n$, the degree of the Fourier extension, and $T$, the size of the extension interval, in the case of $\tilde{f}_{n,m}$.  Hence, given $m$, we need some way to select these values in order to make a comparison.  We shall do this as follows.  In the case of $f_{n,m}$ we let $n$ be as large as possible whilst keeping the condition number $\kappa_{\mathrm{PLS}} \leq \kappa_0$, where $\kappa_0 $ is some fixed value (we use $\kappa_ 0 =10$ in our experiments), and for $\tilde{f}_{n,m}$ we first fix $T$ and then proceed in the same way by choosing $n$ such that $\kappa_{\mathrm{FE}} \leq \kappa_0$.  In other words, for both approximations we ensure that the condition number is no worse than $\kappa_0$.  Thus both methods are guaranteed to be equally robust with respect to perturbations.

Note that computing the $\kappa_{\mathrm{PLS}}$ is straightforward: provided an orthonormal polynomial basis of scaled Legendre polynomials is used, one only needs to determine the minimal singular value of a particular matrix \cite{BAACHAccRecov}.  On the other hand, computing the condition number for the Fourier extension is more challenging.  As was explained in \cite{FEStability}, one cannot determine the condition number by simply examining singular values of the matrix of the linear system resulting from \R{FE}.  Instead, we compute the condition number as follows.  Since the method is linear, we have
\bes{
\kappa_{\mathrm{FE}} = \sup_{\substack{b \in \bbC^{2m+1} \\ b \neq 0}} \frac{\| F_{n,m}(b) \|_{2} }{\| b \|_{l^2}},
}
where $F_{n,m}(b) := \underset{\phi \in \cS_{n}}{\operatorname{argmin}} \{ \sum_{|j| \leq m} | b_j - \hat{\phi}_j |^2 \}$ is the Fourier extension computed from the vector of Fourier coefficients $b$.  We can therefore approximate $\kappa_{\mathrm{FE}}$ by randomly selecting vectors $b$ and computing $\| F_{n,m}(b) \|_2 / \| b \|_{l^2}$.  If $b^{[1]},\ldots,b^{[t]} \in \bbC^{2m+1}$ are chosen uniformly at random with $\| b^{[j]} \|_{l^2} \leq 1$, $j=1,\ldots,t$, we consider the approximation
\bes{
 \tilde\kappa_{\mathrm{FE}} = \max_{j=1,\ldots,t } \frac{\| F_{n,m}(b^{[j]}) \|_{2} }{\| b^{[j]} \|_{l^2}} \approx \kappa_{\mathrm{FE}}.
}
In all our experiments we take the number of trials $t=100$.  In order to compute $\tilde \kappa_{\mathrm{FE}}$ we also need to approximate $\| F_{n,m}(b^{[j]}) \|_{2} $.  We do this with an equispaced quadrature based on $2001$ nodes.

\begin{figure}
\begin{center}
$\begin{array}{ccc}
 \includegraphics[width=4.75cm]{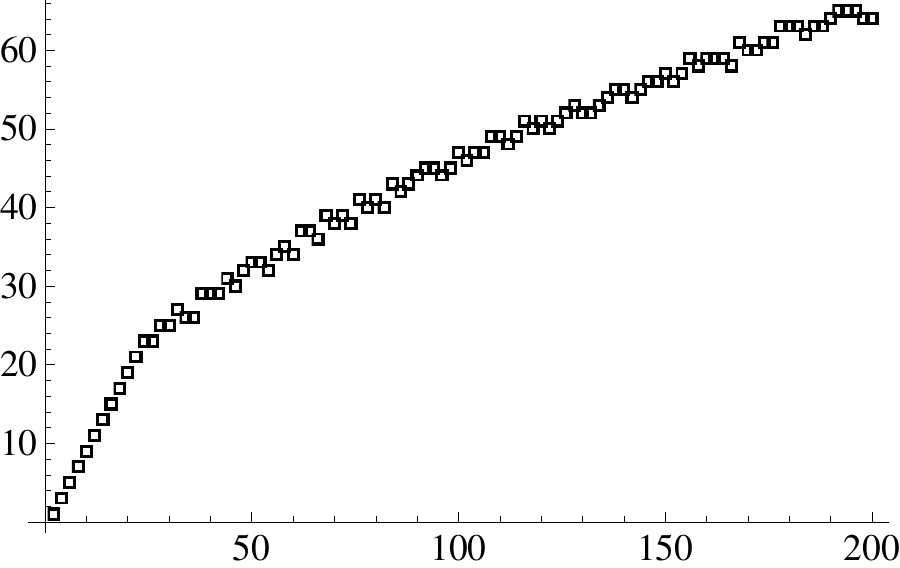} & \hspace{2pc} &  \includegraphics[width=4.75cm]{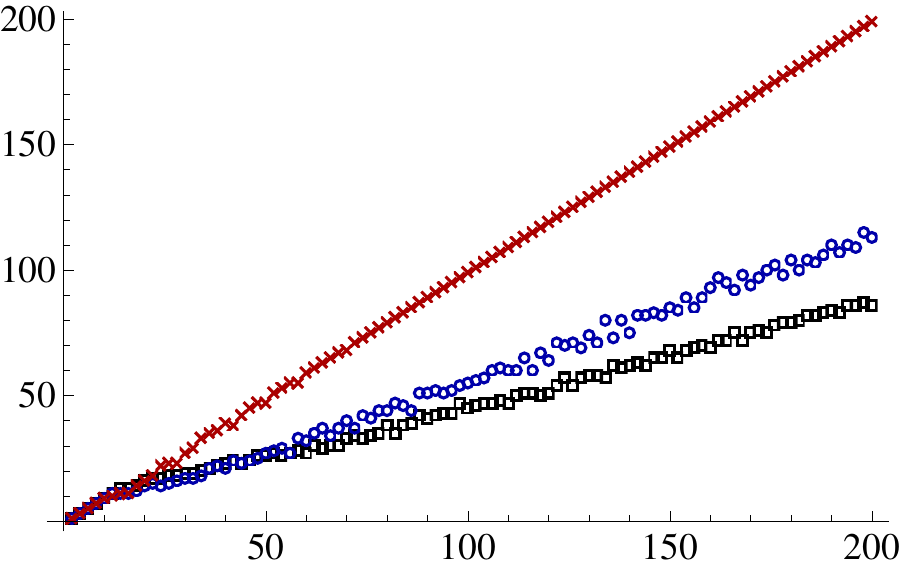} \\
 \includegraphics[width=4.75cm]{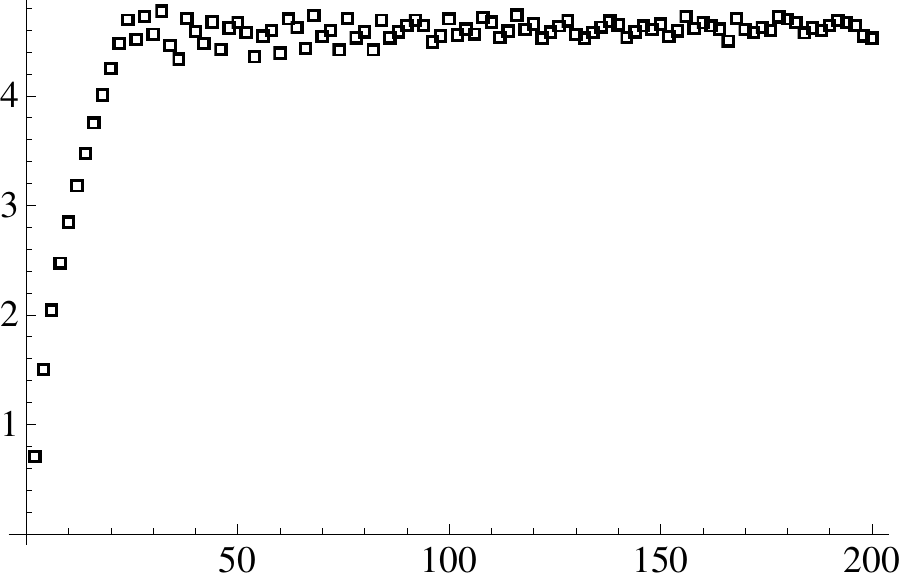} & \hspace{2pc} &  \includegraphics[width=4.75cm]{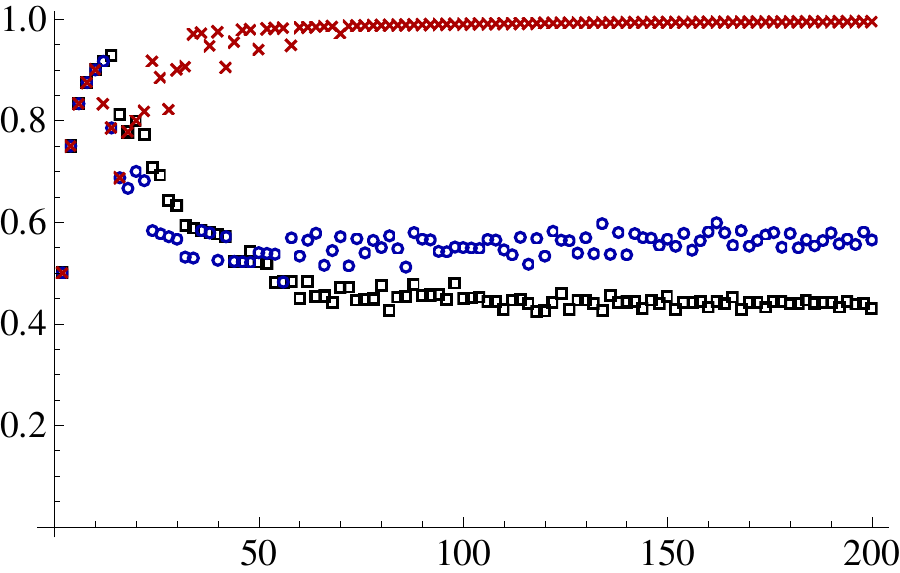} \\
 \mbox{polynomial least squares} && \mbox{Fourier extension}
\end{array}$
\caption{\small Top row: The parameter $n$ against $m=1,\ldots,200$, where $n$ is chosen as large as possible so that $\tilde \kappa_{\mathrm{FE}} \leq 10$ for the Fourier extension and $\kappa_{\mathrm{PLS}} \leq 10$ for polynomial least squares.  In the right column, squares, circles and crosses correspond to the parameter values $T=\frac32,2,4$ respectively.  Bottom row: the scaled parameter $n / \sqrt{m}$ (left) and $n/m$ (right).}\label{f:SSR}
\end{center}
\end{figure}

In Figure \ref{f:SSR} we plot the computed values of $n$ against the number of samples $m$  for each method.  Note that $n = \ord{\sqrt{m}}$ for polynomial least squares, exactly as the results of the last section predict.  On the other hand $n$ grows linearly in $m$ for the Fourier extension, with the constant of proportionality depending on the choice of the parameter $T$.

\begin{figure}
\begin{center}
$\begin{array}{ccc}
 \includegraphics[width=4.75cm]{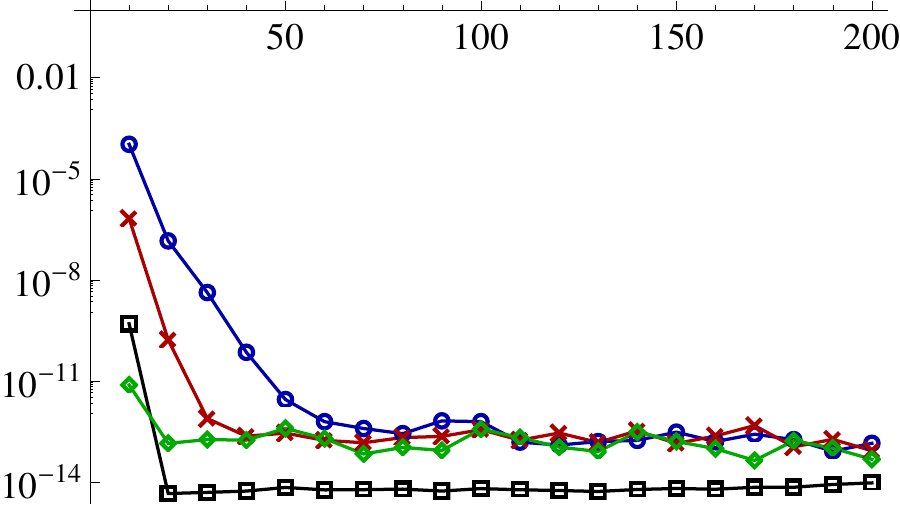} & \hspace{2pc} &  \includegraphics[width=4.75cm]{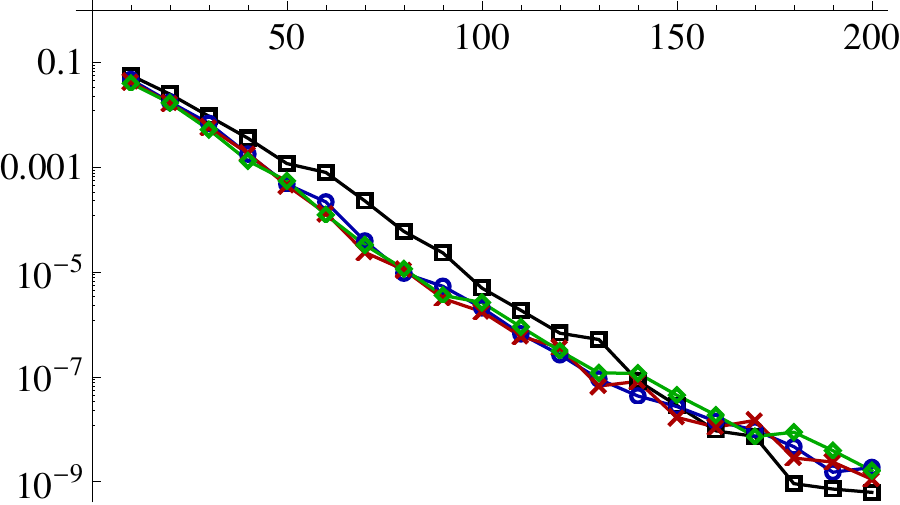} \\
f(x) =  \E^{(x-1)} & &f(x) =  \E^{100 (x-1) } \\\\
 \includegraphics[width=4.75cm]{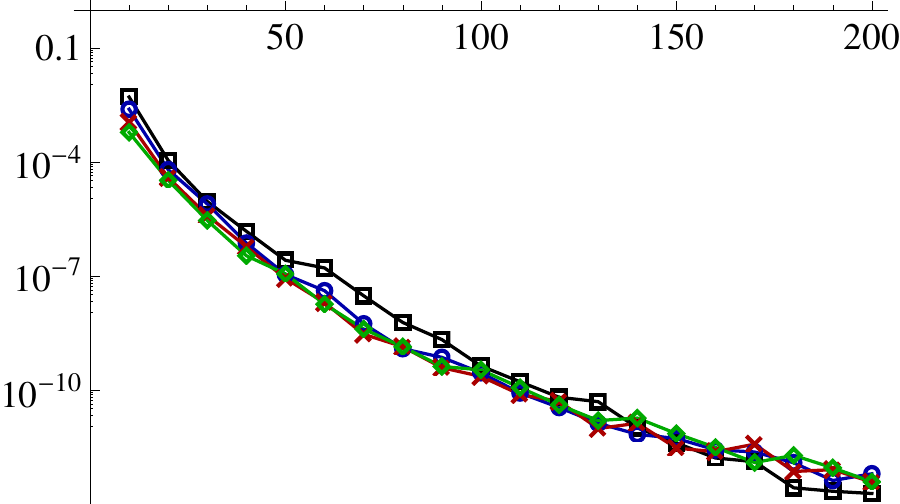} & \hspace{2pc} &  \includegraphics[width=4.75cm]{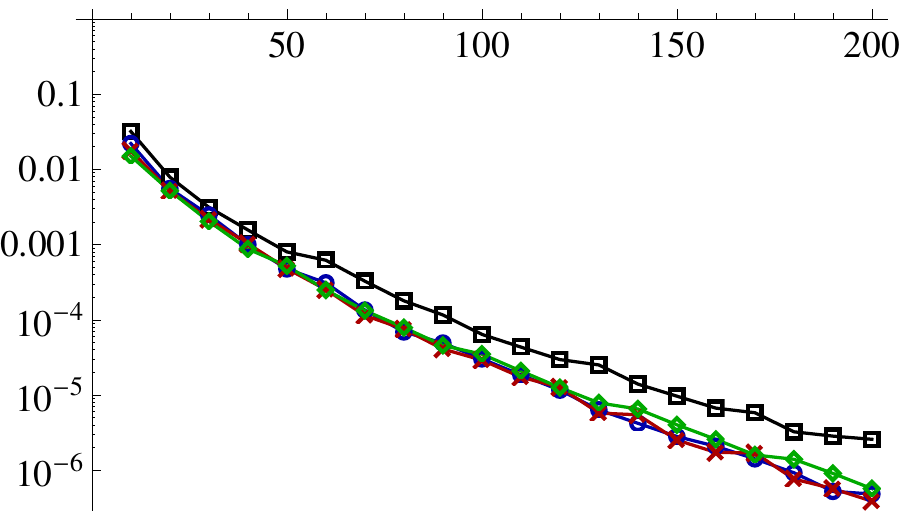}\\
 f(x) = \frac{1}{10-9 x} & & f(x) = \frac{1}{50-49 x}\\\\ \includegraphics[width=4.75cm]{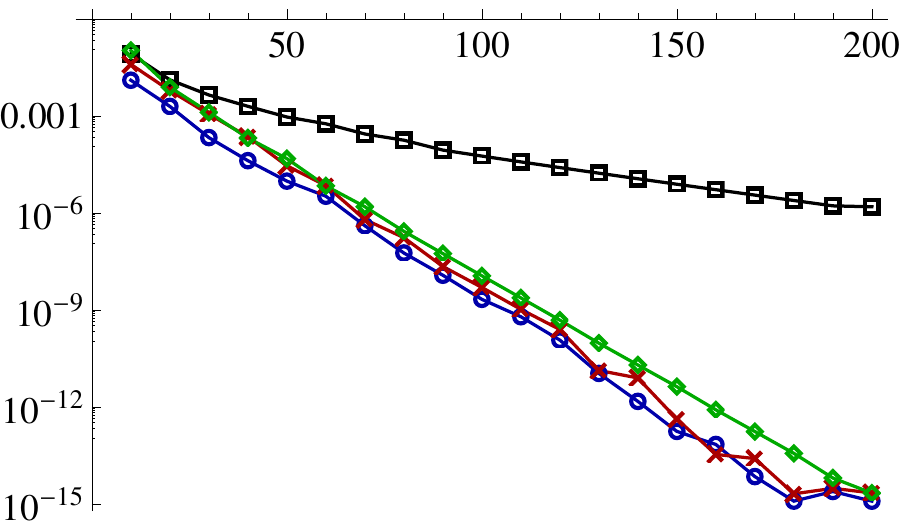} & \hspace{2pc} &  \includegraphics[width=4.75cm]{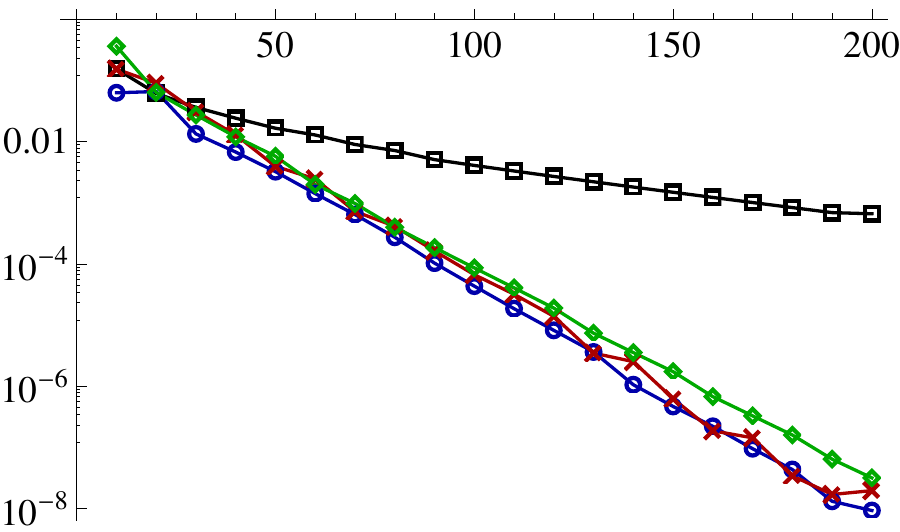} \\
f(x) =  \frac{1}{1+25 x^2} & & f(x) = \frac{1}{1+100 x^2} \\\\
 \includegraphics[width=4.75cm]{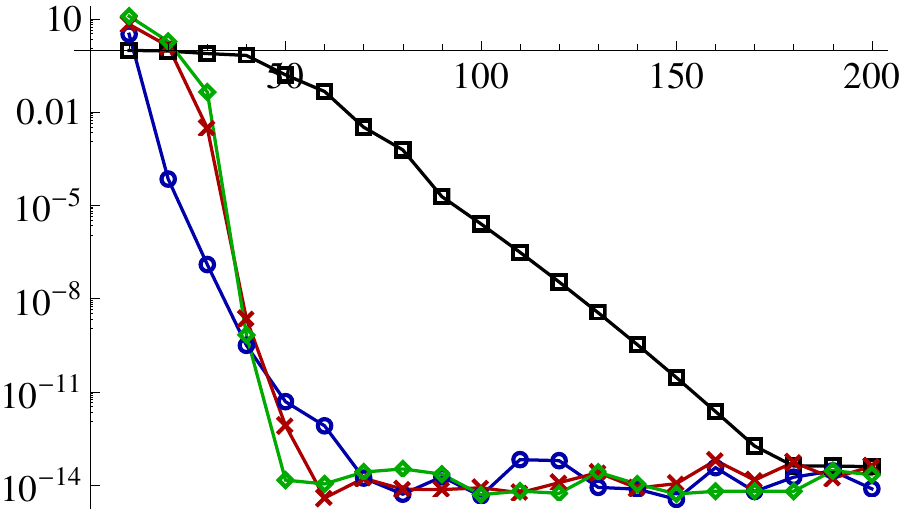} & \hspace{2pc} &  \includegraphics[width=4.75cm]{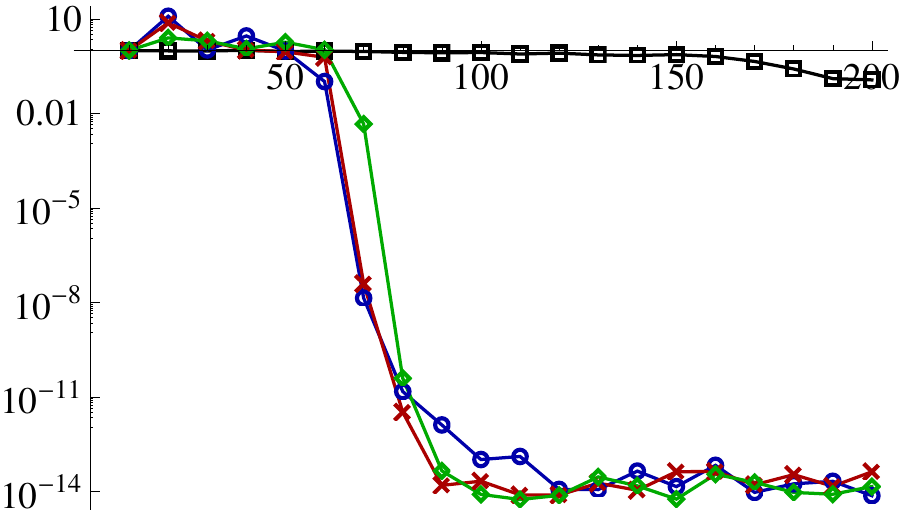}\\
f(x) =  \cos 7 \sqrt{2} \pi x & & f(x) = \cos 14 \sqrt{2} \pi x 
\end{array}$
\caption{\small Errors for the polynomial least squares and Fourier extension approximations $f_{n,m}$ and $\tilde{f}_{n,m}$ against $m=10,20,\ldots,200$.  Squares correspond to polynomial least squares, circles, crosses and diamonds correspond to Fourier extensions with parameter $T = \frac{3}{2},2,4$ respectively.  For each $m$, the parameter $n$ was determined so that the condition number of the corresponding method was at most $10$.}  \label{f:fnapproximation}
\end{center}
\end{figure}

We now use these values to compare the two methods for a variety of different test functions.  In Figure \ref{f:fnapproximation} we consider the following four generic types of functions:
\begin{enumerate}
\item[(i)] Entire functions with boundary layers, e.g.\ $f(x) = \E^{a (x-1)}$, $a \gg 1$.
\item[(ii)] Meromorphic functions with real singularities near $x=1$, e.g.\ $f(x) = \frac{1}{a+1-a x}$, $a \gg 1$.
\item[(iii)] Meromorphic functions with complex singularities on the imaginary axis near $x=0$, e.g.\ the classical Runge example $f(x) = \frac{1}{1+a^2 x^2}$, $a \gg 1$.
\item[(iv)] Entire, oscillatory functions, e.g.\ $f(x) = \cos \omega \pi x$, $\omega \gg 1$.
\end{enumerate}

Note that algebraic polynomials are particularly well suited to (i) and (ii).  One can resolve a boundary layer of width $1/a$ using a polynomial of degree $\ord{\sqrt{a}}$ \cite{boyd,naspec}.  Fourier approximations such as Fourier extensions are less well suited for boundary layers, since they require $\ord{a}$ degrees of freedom.  However, when only Fourier data is prescribed the polynomial least squares method can stably reconstruct a polynomial of degree at most $\ord{\sqrt{m}}$, whereas the Fourier extension method recovers a Fourier series of degree $\ord{m}$.  This means that a boundary layer of width $a$ cannot be resolved any more efficiently by polynomial least squares than by the Fourier extension when applied to this type of data.  Hence, this particular advantage of using a polynomial reconstruction space disappears.  As can be seen in Figure \ref{f:fnapproximation}, $f_{n,m}$ and $\tilde{f}_{n,m}$ give roughly the same errors for functions of type (i) and (ii).

This figure also shows that Fourier extensions significantly outperform polynomial least squares for functions of type (iii) and (iv).  This can be explained in a similar manner.  Up to constant factors, Fourier extensions and algebraic polynomial approximations require the same number of degrees of freedom $n$ to resolve functions of type (iii) and (iv) \cite{BADHFEResolution}.  However, since we recover a Fourier extension of degree $\ord{m}$, whereas we can only recover a polynomial of $\ord{\sqrt{m}}$, we obtain a vastly superior approximation using the former.

Notice also one interesting facet of Figure \ref{f:fnapproximation}: namely, the choice of $T$ has little effect on the approximation error.  This is due to the balance between degrees of freedom (i.e.\ $n$) and approximation properties.  Larger $T$ means a larger $n$ can be used for a given $m$ whilst retaining the condition number (see Figure \ref{f:SSR}).  However, larger $T$ also translates into slower (although still superalgebraic/geometric) convergence in $n$ \cite{BADHFEResolution}.  Perhaps surprisingly, these two effects balance in practice, rendering the choice of $T$ insignificant.  To the best of our knowledge, this observation has not previously been made in the literature on Fourier extensions.

To summarize, Fourier extensions appear to present an attractive method for the problem of recovering analytic functions from Fourier data.  In particular, they allow one to circumvent the stability barrier of Theorem \ref{t:ctsPTK} to a substantial extent, since they only converge down to a finite accuracy (albeit on the order of machine precision).  A more thorough comparison, incorporating more of the methods discussed in Section \ref{s:methods}, is a topic for future work.

\section*{Acknowledgements}
The authors would like to thank Karlheinz Gr\"ochenig, Tomasz Hrycak, Daan Huybrechs,  Nilima Nigam, Rodrigo Platte and Nick Trefethen for helpful discussions and comments.

\bibliographystyle{abbrv}
\small
\bibliography{GibbsStabilityRefs}

\end{document}